\documentclass[10pt]{amsart}
\usepackage{amsmath}
\usepackage{amsfonts}
\usepackage{amsthm}
\usepackage{enumerate}
\usepackage{amssymb}
\usepackage{amscd}

\DeclareMathOperator{\gr}{gr}

\DeclareMathOperator{\gld}{gld}
\DeclareMathOperator{\Aut}{Aut}

\DeclareMathOperator{\Ext}{Ext}
\DeclareMathOperator{\pd}{pd}

\DeclareMathOperator{\ad}{ad}
\DeclareMathOperator{\Cdim}{Cdim}
\DeclareMathOperator{\GK}{GK}
\DeclareMathOperator{\Ann}{Ann}
\DeclareMathOperator{\Ch}{Ch}

\begin{document}
\theoremstyle{plain}
\newtheorem{quest}{Question}
\renewcommand{\thequest}{\Alph{quest}}
\newtheorem*{trm}{Theorem}
\newtheorem*{lem}{Lemma}
\newtheorem*{prop}{Proposition}
\newtheorem*{conj}{Conjecture}
\newtheorem*{rems}{Remarks}
\newtheorem*{thm}{Theorem}
\newtheorem*{example}{Example}
\newtheorem*{examples}{Examples}
\newtheorem*{cor}{Corollary}
\newtheorem*{hyp}{Hypothesis}
\newtheorem*{thrm}{Theorem}
\theoremstyle{remark}
\newtheorem{defn}{Definition}
\newtheorem*{rem}{Remark}
\newtheorem*{notn}{Notation}
\newcommand{\Fp}{\mathbb{F}_p}
\newcommand{\Zp}{\mathbb{Z}_p}
\newcommand{\Qp}{\mathbb{Q}_p}
\newcommand{\Kr}{\mathcal{K}}
\newcommand{\Rees}[1]{\widetilde{#1}}
\newcommand{\invlim}{\lim\limits_{\longleftarrow}}

\title[Iwasawa algebras survey]{Ring-theoretic properties \\ of Iwasawa algebras: a survey}

\author{K. Ardakov and K.A. Brown }

\address{Ardakov:  DPMMS, University of Cambridge,
Centre for Mathematical Sciences, Wilberforce Road, Cambridge CB3 0WB, UK}

\email{K.Ardakov@dpmms.cam.ac.uk}

\address{Brown: Department of Mathematics,
University of Glasgow, Glasgow G12 8QW, UK}

\email{kab@maths.gla.ac.uk}

\begin{abstract}
This is a survey of the known properties of Iwasawa algebras,
which are completed group rings of compact $p-$adic analytic
groups with coefficients the ring $\Zp$ of $p-$adic integers or
the field $\Fp$ of $p$ elements. A number of open questions are
also stated.
\end{abstract}



\thanks{ Some of the work for this article was done in June
2005, when Ardakov was visiting the University of Glasgow with the support of the Edinburgh Mathematical Society
Research Support Fund and the Glasgow Mathematical Journal Learning and Research Support Fund.}

\maketitle
\section{Introduction}

Noncommutative Iwasawa algebras form a large and interesting class
of complete semilocal noetherian algebras, constructed as
completed group algebras of compact $p-$adic analytic groups. They
were defined and their fundamental properties were derived in M.
Lazard's monumental 1965 paper \cite{L}, but in the twenty years
from 1970 they were little studied. Interest in them has been
revived by developments in number theory over the past fifteen
years, see for example \cite{C}. Prompted by this renewed
interest, and helped of course by the better understanding of
noncommutative noetherian algebra gained since 1965, a number of
recent papers have built on Lazard's initial work. The emerging
picture is of a class of rings which in some ways look similar to
the classical commutative Iwasawa algebras, (which are rings of
formal power series in finitely many commuting variables over the
$p-$adic integers), but which in other respects are very different
from their commutative counterparts. And while some progress has
been made in understanding these rings, many aspects of their
structure and representation theory remain mysterious.

It is the purpose of this article to provide a report of what is
known about Iwasawa algebras at the present time, and to make some
tentative suggestions for future research directions. We approach
the latter objective through the listing of a series of open
questions, scattered throughout the text. In an attempt to make
the paper accessible to readers from as wide a range of
backgrounds as possible, we have tried to give fairly complete
definitions of all terminology; on the other hand, most proofs are
omitted, although we have tried to give some short indication of
their key points where possible. An exception to the omission of
proofs occurs in the discussion of maximal orders in
(\ref{max})-(\ref{Suff}), where we include some original material.

Fundamental definitions and examples are given in Section 2; in
particular we recall the definition of a \emph{uniform} pro$-p$
group in (\ref{uniform}), and make the important observation
(\ref{iwa})(1) that every Iwasawa algebra can be viewed as a
crossed product of the Iwasawa algebra of a uniform group by a
finite group. This has the effect of focusing attention on the
Iwasawa algebra of a uniform group - this is filtered by the
powers of its Jacobson radical, and the associated graded algebra
is a (commutative) polynomial algebra. This fact and its
consequences for the structure of the Iwasawa algebras of uniform
groups are explored in Section 3; then in Section 4 we examine how
properties of general Iwasawa algebras can be deduced from the
uniform case using (\ref{iwa})(1). Section 5 concerns dimensions:
first, the global (projective) dimension and the injective
dimension, whose importance is enhanced because Iwasawa algebras
satisfy the \emph{Auslander-Gorenstein condition}, whose
definition and properties we recall. In particular,
Auslander-Gorenstein rings possess a so-called \emph{canonical
dimension function}; we explain this and describe some of the
properties of the canonical dimension of an Iwasawa algebra in
(\ref{dim})-(\ref{char}). The Krull-Gabriel-Rentschler dimension is
discussed in (\ref{KdimOmega}). Finally, our very sparse knowledge
of the two-sided ideals of Iwasawa algebras is summarised in
Section 6.

\section{Key definitions} \label{defns} Iwasawa algebras are completed group algebras. We begin by recalling which
groups are involved, then give the definition of the algebras.

\subsection{Compact $p-$adic analytic groups} \label{groups} Let $p$ be a prime integer and let $\Zp$ denote
the ring of $p-$adic integers. A group $G$ is {\it compact $p-$adic analytic} if it is a topological group which
has the structure of a $p-$adic analytic manifold - that is, it has an atlas of open subsets of $\Zp^{n}$, for
some $n \geq 0.$ Such groups can be characterised in a more intrinsic way, thanks to theorems due to Lazard,
dating from his seminal 1965 paper \cite{L}. Namely, a topological group $G$ is compact $p-$adic analytic if and
only if $G$ is profinite, with an open subgroup which is pro$-p$ of finite rank, if and only if $G$ is a closed
subgroup of $GL_d(\Zp)$ for some $d \geq 1.$ Nowadays, these equivalences are usually viewed as being consequences
of deep properties of finite $p-$groups; a detailed account from this perspective can be found in \cite[Part
II]{DDMS}.

\bigskip

\noindent \textbf{Examples:} (1) Every finite group is $p-$adic analytic, for every prime $p.$

(2) The abelian $p-$adic analytic groups are the direct products of finitely many copies of the additive group of
$\Zp$ with a finite abelian group \cite[page 36]{DDMS}.

(3) For any positive integer $d$ the groups $GL_d(\Zp)$ and $SL_d(\Zp)$ are compact $p-$adic analytic. More
generally, given any root system $X_{\ell}$ one can form the {\it universal Chevalley group} $\mathcal{G}_{\Zp}(X_{\ell})$,
\cite[page 353]{DDMS}. This is a compact $p-$adic analytic group. For more information about Chevalley groups, see
\cite{Carter}.

(4) Let $d$ and $t$ be positive integers. The {\it $t$-th congruence subgroup in} $SL_d(\Zp)$ is the kernel
$\Gamma_t(SL_d(\Zp))$ of the canonical epimorphism from $SL_d(\Zp)$ to $SL_d(\Zp/p^t\Zp).$ One sees at once from
the equivalences above that $\Gamma_t(SL_d(\Zp))$ is compact $p-$adic analytic, as indeed are
$\Gamma_t(GL_d(\Zp))$ and $\Gamma_t(\mathcal{G}_{\Zp}(X_{\ell}))$ for any root system $X_{\ell}$.

\bigskip

\noindent \textbf{Notation:} When discussing a topological group $G$ we shall use $\overline{H}$ to denote the
closure of a subset $H$ of $G$ in $G$; and when we refer to, say, $G$ as being {\it generated by} elements $\{g_1,
\ldots , g_d \}$ we mean that $G = \overline{ \langle g_1 , \ldots , g_d \rangle}.$ In particular, $G$ is {\it
finitely generated} if $G = \overline{\langle X \rangle}$ for a finite subset $X$ of $G.$ For a subset $X$ of $G$,
$X^p$ denotes the subgroup of $G$ generated by the subset $\{x^p : x \in X \}$ of $G$.

\subsection{Iwasawa algebras}\label{iwa} Let $G$ be a compact $p-$adic
analytic group. The {\it Iwasawa algebra of $G$} is
$$ \Lambda_G \quad := \quad \invlim \Zp [G/N], $$
where the inverse limit is taken over the open normal subgroups $N$ of $G$. Closely related to $\Lambda_G$ is its
epimorphic image $\Omega_G$, defined as
$$  \Omega_G \quad = \quad \invlim \Fp [G/N], $$
where $\Fp$ is the field of $p$ elements. Often, a property of $\Lambda_G$ can easily be deduced from the
corresponding property of $\Omega_G$, and vice versa; where this is routine we will frequently save space by
stating only one of the two variants.

Suppose that $H$ is an open normal subgroup of $G.$ Let $\mathcal{C}_H$ denote the set of open normal subgroups of
$G$ which are contained in $H$. Then clearly $\Lambda_{G} = {\lim\limits_{\longleftarrow}} \Zp [G/U]$ where $U$
runs over $\mathcal{C}_H$, and so it follows at once that
\begin{equation} \label{crossed} \Lambda_G \quad \cong \quad \Lambda_H \ast (G/H),  \tag{1}
\end{equation}
a crossed product of $\Lambda_H$ by the finite group $G/H.$ We shall see that, combined with a judicious choice of
the subgroup $H$, the isomorphism (\ref{crossed}) reduces many questions about $\Lambda_G$ and $\Omega_G$ to the
analysis of certain crossed products of finite groups. Usually, the right subgroup $H$ to choose is a {\it
uniform} one, defined as follows.

\subsection{Uniform groups} \label{uniform}  Let $G$ be a pro$-p$ group.
Define $P_1(G) = G$ and $P_{i+1}(G) =
\overline{P_i(G)^p[P_i(G),G]}$ for $i \geq 1.$ The decreasing
chain of characteristic subgroups $$ G = P_1(G) \supseteq P_2(G)
\supseteq \cdots \supseteq P_i(G) \supseteq \cdots \supseteq
\cap_{i=1}^{\infty}P_i(G) = 1 $$ is called the {\it lower
$p-$series} of $G.$ The group $G$ is \emph{powerful} if
$G/\overline{G^p}$ is abelian (for $p$ odd), or $G/\overline{G^4}$
is abelian (when $p = 2).$ Finally, $G$ is {\it uniform} if it is
powerful, finitely generated, and $$ |G : P_2(G)| = |P_i(G) :
P_{i+1}(G)| $$ for all $i \geq 1.$

Now we can add one further characterisation, also essentially due to Lazard, to those given in (\ref{groups}): a
topological group $G$ is compact $p-$adic analytic if and only if it has an open normal uniform pro$-p$ subgroup
of finite index, \cite[Corollary 8.34]{DDMS}.

\bigskip

\noindent \textbf{Examples:} (1) Of course, $(\Zp)^{\oplus d}$ is uniform for all $d \geq 1.$

(2)The groups $\Gamma_1(GL_d(\Zp))$ (for $p$ odd) and $\Gamma_2(GL_d(\mathbb{Z}_2))$ are uniform \cite[Theorem 5.2]{DDMS}.

\bigskip

Let $G$ be uniform, with $|G : P_2(G)| = p^d.$ The non-negative integer $d$ is called the \emph{dimension} of $G$;
it is equal to the cardinality of a minimal set of (topological) generators of $G$, \cite[Definition 4.7 and
Theorem 3.6]{DDMS}. More generally, we can define the dimension of an arbitrary compact $p-$adic analytic 
group to be the dimension of any open uniform subgroup; this is unambiguous \cite[Lemma 4.6]{DDMS}, and coincides
with the dimension of $G$ as a $p-$adic analytic manifold, \cite[Definition 8.6 and Theorem 8.36]{DDMS}.

\subsection{Completed group algebras}\label{comp} In fact $\Lambda_G$
and $\Omega_G$ are $I-$adic completions of the ordinary group
algebras $\Zp[G]$ and $\Fp [G]$, for suitable choices of ideals
$I$. It is most convenient for us to state the result for uniform
groups, although it can obviously be extended to the general case
using (\ref{iwa})(\ref{crossed}).

\begin{thm} Let $G$ be a uniform pro$-p$ group, and let I denote
the augmentation ideal of $\Fp [G].$ Then $\Omega_G$ is isomorphic
to the $I-$adic completion of $\Fp [G].$ There is a similar result
for $\Zp [G].$
\end{thm}

Indeed the theorem follows quite easily from the observations that
the lower $p-$series $P_i(G)$ is coterminal with the
family of all open normal subgroups of $G$, and that the powers of
$I$ are coterminal with the ideals of $\Fp [G]$ generated by the
augmentation ideals of the subgroups $P_i(G)$, \cite[\S
7.1]{DDMS}.

\section{The case when $G$ is uniform}

Throughout this section, we assume that $G$ is a uniform pro-$p$ group of dimension $d$.
We fix a topological generating set $\{a_1,\ldots,a_d\}$ for $G$.

\subsection{The ``PBW" Theorem}\label{PBW}

It follows at once from Theorem (\ref{comp}) that the usual group
algebra $\Fp[G]$ embeds into $\Omega_G$. For $i=1,\ldots,d$, let
$b_i = a_i-1 \in \Fp[G] \subseteq \Omega_G$. Then we can form
various monomials in the $b_i$: if $\alpha =
(\alpha_1,\ldots,\alpha_d)$ is a $d-$tuple of nonnegative
integers, we define
\[\mathbf{b}^\alpha = b_1^{\alpha_1}\cdots b_d^{\alpha_d} \in \Omega_G.\]
Note that this depends on our choice of ordering of the $b_i$'s, because $\Omega_G$ is noncommutative unless $G$
is abelian. The following basic result shows that $\Omega_G$ is a ``noncommutative formal power series ring"; it
follows from the strong constraints which the hypothesis of uniformity imposes on the quotients
$P_i(G)/P_{i+1}(G)$ of $G,$ \cite[Theorem 7.23]{DDMS}.

\begin{thm}

Every element $c$ of $\Omega_G$ is equal to the sum of a uniquely determined convergent series
\[c = \sum_{\alpha \in\mathbb{N}^d} c_{\alpha} \mathbf{b}^{\alpha} \]
where $c_\alpha \in \Fp$ for all $\alpha \in \mathbb{N}^d$.
\end{thm}

We record an immediate consequence of both this result and of Theorem (\ref{comp}):
\begin{cor} The Jacobson radical $J$ of $\Omega_G$ is equal to
\[J = b_1\Omega_G + \cdots + b_d \Omega_G = \Omega_G b_1 + \cdots + \Omega_G b_d.\]
Hence $\Omega_G$ is a scalar local ring with $\Omega_G/J \cong \Fp$.
\end{cor}
\begin{proof}
If $c\in \Omega_G$ is such that $c_0\neq 0$, then $1 - c$ is invertible with inverse $1 + c + c^2 + \cdots \in \Omega_G$.
\end{proof}

Theorem (\ref{PBW}) says that the monomials $\{\mathbf{b}^\alpha : \alpha \in \mathbb{N}^d\}$ form a topological
basis for $\Omega_G$, and is thus analogous to the classical Poincar\'e-Birkhoff-Witt theorem for Lie
algebras $\mathfrak{g}$ over a field $k$ which gives a vector space basis for the universal enveloping algebra
$\mathcal{U}(\mathfrak{g})$ in terms of monomials in a fixed basis for $\mathfrak{g}$ \cite{Dixm}. Nevertheless we
should bear in mind that explicit computations in $\Omega_G$ are often much more difficult than those in
$\mathcal{U}(\mathfrak{g})$, since the Lie bracket of two generators $b_i$, $b_j$ for $\Omega_G$ is in general an
infinite power series with obscure coefficients.

\subsection{Example}
\label{SL2} Let $p$ be odd for simplicity and let $G = \Gamma_1
(SL_2(\Zp))$ be the first congruence kernel of $SL_2(\Zp)$. Then
\[a_1 = \begin{pmatrix} \exp(p) & 0 \\ 0 & \exp(-p)\end{pmatrix}, \quad a_2 = \begin{pmatrix} 1 & p \\ 0 & 1\end{pmatrix},\quad a_3 = \begin{pmatrix} 1 & 0 \\ p & 1\end{pmatrix}.\]
is a topological generating set for $G$. Setting $b_i = a_i - 1$, elementary (but tedious) computations yield
\[
\begin{array}{llll}
[b_1,b_2] &\equiv &2b_2^p &\mod J^{p+1}\\

[b_1,b_3] &\equiv &-2b_3^p &\mod J^{p+1}\\

[b_2,b_3] &\equiv &b_1^p &\mod J^{p+1}.
\end{array}
\]
Here $J = b_1\Omega_G + b_2\Omega_G + b_3\Omega_G$ denotes the Jacobson radical of $\Omega_G$. Using Proposition (\ref{skew}) it is possible to produce more terms in the power series expansion of $[b_1,b_2]$ and $[b_1,b_3]$. However, we consider $[b_2,b_3]$ to be inaccessible to computation.

\subsection{Skew power series rings}
\label{skew} It is well known that if $\mathfrak{g}$ is a finite dimensional soluble Lie algebra over a field
$k$, then its universal enveloping algebra $\mathcal{U}(\mathfrak{g})$ can be thought of as an ``iterated skew
polynomial ring":
\[\mathcal{U}(\mathfrak{g}) \cong k[x_1;\sigma_1,\delta_1][x_2;\sigma_2,\delta_2]\cdots[x_n;\sigma_n,\delta_n]\]
for some appropriate automorphisms $\sigma_i$ and derivations $\delta_i$ (in fact, the $\sigma_i$s can be chosen
to be trivial). This is because any such Lie algebra $\mathfrak{g}$ has a chain of subalgebras
\[0 = \mathfrak{h}_0 \subset \mathfrak{h}_1 \subset \mathfrak{h}_2 \subset \cdots \subset \mathfrak{h}_n = \mathfrak{g}\]
with $\mathfrak{h}_{i-1}$ an ideal in $\mathfrak{h}_i$, so choosing some $x_i \in\mathfrak{h}_i \backslash
\mathfrak{h}_{i-1}$ ensures that
\[\mathcal{U}(\mathfrak{h}_i) \cong \mathcal{U}(\mathfrak{h}_{i-1})[x_i;\delta_i]\]
where $\delta_i$ is the derivation on $\mathcal{U}(\mathfrak{h}_{i-1})$ extending $\ad(x_i)_{\mid
\mathfrak{h}_{i-1}}$.

An analogous result holds for Iwasawa algebras. More precisely, we have the

\begin{prop}
Suppose that $G$ has closed normal subgroup $H$ such that $G/H \cong \Zp$. Then $\Omega_G$ is a skew power series ring with coefficients in $\Omega_H$:
\[\Omega_G \cong \Omega_H[[t;\sigma,\delta]].\]
\end{prop}
\begin{proof} See \cite[\S 4]{VS}.
\end{proof}
Schneider and Venjakob \cite{VS} establish a general theory of skew power series rings $S = R[[t; \sigma,\delta]]$
over a pseudocompact ring $R$. Here $\sigma$ can be any topological automorphism of $R$ and $\delta$ is a
$\sigma-$derivation in the sense of \cite[1.2.1]{MR}, satisfying some extra conditions which are required to make
the relation
\[ ta = \sigma(a)t + \delta(a)\]
extend to a well-defined multiplication on $S$.

Consequently, the Iwasawa algebra $\Omega_G$ of any soluble uniform pro-$p$ group $G$ can
be thought of as an iterated skew power series ring over $\Fp$.

For example, in Example \ref{SL2}, the topological subring of $\Omega_G$ generated by $b_1$ and $b_2$ is actually
the Iwasawa algebra $\Omega_B$ where $B = \overline{\langle a_1,a_2\rangle}$ is a Borel subgroup of $G$. Since $B$
is soluble with closed normal subgroup $\overline{\langle a_2\rangle}$, $\Omega_B$ is isomorphic to the skew power
series ring $\Fp[[b_2]][[b_1;\sigma,\delta]]$ for some appropriate $\sigma$ and $\delta$. This justifies the claim
that the commutator of $b_1$ and $b_2$ is at least partially accessible to computation.

There is surely considerable scope to develop further the ``abstract'' theory of skew power series algebras
initiated in \cite{VS} - for instance, one could easily pose skew power series versions of a number of the
questions we list later, in Section 6. As a prompt for more work, here are two ``general'' questions:

\begin{quest}(1) Are there conditions on $R, \sigma$ and $\delta$ such that $S = R[[t; \sigma,\delta]]$ can be
described without involving a derivation - that is, as $S = R'[[t';\sigma']]$, possibly after some Ore
localisation?\footnote{Compare with \cite{Cau}.}

\noindent (2) Are there conditions on $R, \sigma$ and $\delta$ such that every two-sided ideal of the skew power
series ring $S = R[[t; \sigma,\delta]]$ is generated by central elements and ``polynomial'' elements\footnote{By
the latter, we mean elements of $R[t; \sigma,\delta].$ }?
\end{quest}

\subsection{The $J-$adic filtration}
We remind the reader that a \emph{filtration} on a ring $R$ is an ascending sequence
\[\cdots \subseteq F_iR \subseteq F_{i+1}R \subseteq \cdots\]
of additive subgroups such that $1 \in F_0R$, $F_iR.F_jR \subseteq F_{i+j}R$ for all $i,j \in \mathbb{Z}$, and $\cup_{i\in\mathbb{Z}} F_iR = R$.

Let $J$ denote the Jacobson radical of $\Omega_G$. The $J-$adic
filtration on $\Omega_G$ is defined as follows: $F_i\Omega_G =
J^{-i}$ for $i\leq 0$ and $F_i\Omega_G = \Omega_G$ f or $i\geq 0$;
this is an example of a \emph{negative} filtration. The basic tool
which allows one to deduce many ring-theoretic properties of
Iwasawa algebras is the following result, which can be deduced
from Theorem (\ref{PBW}), see \cite[Theorem 7.24 and remarks on
page 160]{DDMS}. We denote the associated graded ring $
\bigoplus_{i \in \mathbb{Z}}F_{i+1}\Omega_G/F_i\Omega_G$ by $\gr_J
\Omega_G.$

\begin{thm}
\label{GrR} The graded ring of $\Omega_G$ with respect to the $J-$adic filtration is isomorphic to a polynomial
ring in $d = \dim G$ variables:
\[\gr_J \Omega_G \cong \Fp[X_1,\ldots,X_d].\]
Moreover, $\Omega_G$ is complete with respect to this filtration.
\end{thm}

The $J-$adic filtration is quite different from the filtrations
encountered when studying algebras like universal enveloping
algebras and Weyl algebras, which are nearly always
\emph{positive} (that is, $F_{-1}R = 0$) and often satisfy the
finiteness condition $\dim_k F_iR < \infty$ for all $i\in
\mathbb{Z}$. In particular, there is no well-behaved notion of the
Gel'fand-Kirillov dimension for Iwasawa algebras, a theme we will
return to in \S \ref{Dims}.

However, we are still able to lift many properties of the graded ring back to $\Omega_G$, because the $J-$adic
filtration is {\it complete}, meaning that Cauchy sequences of elements in $\Omega_G$ converge to unique limits.
More precisely, recall \cite[page 83]{LVO} that a filtration on a ring $R$ is said to be \emph{Zariskian},
whenever
\begin{itemize}
\item The Jacobson radical of $F_0R$ contains $F_{-1}R$, and
\item The Rees ring $\widetilde{R} := \bigoplus_{i \in \mathbb{Z}} F_iR\cdot t^i\subseteq R[t,t^{-1}]$ is noetherian.
\end{itemize}

Many filtrations are Zariskian. For example, by \cite[Chapter II, Proposition 2.2.1]{LVO}, any complete filtration
whose associated graded ring is noetherian is necessarily Zariskian. Since any positive filtration is
complete, it follows that if a filtration is positive and has noetherian associated graded ring, then it is
Zariskian. More importantly for us, for any uniform pro$-p$ group $G$, the $J-$adic filtration on $\Omega_G$ is
clearly complete, thanks to Theorem (\ref{comp}); and $\gr_J \Omega_G$ is noetherian by Theorem (\ref{GrR}) and
Hilbert's basis theorem, so the $J-$adic filtration is Zariskian.

\subsection{Lifting information from the graded ring}
\label{lift} We recall here some standard properties of a ring $R$. First, we say that $R$ is \emph{prime} if the
product of any two non-zero ideals of $R$ is again non-zero. By Goldie's theorem \cite[Theorem 2.3.6]{MR}, if $R$
is prime and (right) noetherian then it has a simple artinian classical (right) quotient ring $Q(R).$ If $S$ is
another ring with classical right quotient ring $Q(R),$ so that $Q(R) = Q(S),$ we say that $R$ and $S$ are
\emph{equivalent} if there are units $a,b,c$ and $d$ in $Q(R)$ such that $aRb \subseteq S$ and $cSd \subseteq R.$
Now $R$ is a \emph{maximal (right) order} if it is maximal (with respect to inclusion) within its equivalence
class, \cite[5.1.1]{MR}. (The adjective right is omitted if $R$ is both a maximal right order and a maximal left
order.) The commutative noetherian maximal orders are just the noetherian integrally closed domains \cite[Lemma
5.3.3]{MR}.

The \emph{Krull dimension} $\Kr (M)$ of a finitely generated
(right) module $M$ over a noetherian ring $R$ is a well-defined
ordinal, bounded above by $\Kr (R_R)$; the precise definition can
be found at \cite[6.2.2]{MR}. This concept generalises the
classical commutative definition; like it, it measures the ``size''
of a module and is 0 if and only if the module is non-zero and
artinian.

The \emph{(right) global dimension} of $R$ is defined to be the supremum of the projective dimensions (denoted
$\pd (-)$) of the right $R-$modules, \cite[7.1.8]{MR}. When $R$ is noetherian, its right and left global
dimensions are always equal, \cite[7.1.11]{MR}. We say that $R$ has \emph{finite (right) injective dimension} $d$
if there is an injective resolution of $R_R$ of length $d$, but none shorter. If $R$ is noetherian and has finite
right and left injective dimensions, then these numbers are equal by \cite[Lemma A]{Zaks}. It is also well known \cite[Remark 6.4]{Ven1} that if the (right) global dimension of the noetherian ring $R$ is finite, then it equals the (right) injective dimension of $R$.

It has become apparent over the past 40 years that, when $R$ is noncommutative and noetherian, finiteness of the
injective dimension of $R$ is a much less stringent condition than is the case for commutative noetherian rings -
the structure of (commutative) Gorenstein rings is rich and beautiful. An additional hypothesis which, when
coupled with finite injective dimension, has proved very useful in the noncommutative world is the
\emph{Auslander-Gorenstein condition}. To recall the definition, note first that, for every left $R-$module $M$
and every non-negative integer $i$, $\Ext ^i (M,R)$ is a right $R-$module through the right action on $R.$ The
Auslander-Gorenstein condition on a noetherian ring $R$ requires that, when $M$ is a finitely generated left
$R-$module, $i$ is a non-negative integer and $N$ is a finitely generated submodule of $\Ext^i(M,R),$ then $\Ext^j
(N,R)$ is zero for all $j$ strictly less than $i$; and similarly with ``right" and ``left" interchanged. We say
that $R$ is \emph{Auslander-Gorenstein} if it is noetherian, has finite right and left injective dimensions, and
satisfies the Auslander condition. Commutative noetherian rings of finite injective dimension are
Auslander-Gorenstein. When $R$ is noetherian of finite global dimension and satisfies the Auslander-Gorenstein
condition it is called \emph{Auslander-regular}.

\begin{thm} Let $R$ be a ring endowed with a Zariskian filtration $FR$; then $R$ is necessarily noetherian. Also, $R$ inherits the following properties from $\gr R$:
\begin{enumerate}
\item being a domain,
\item being prime,
\item being a maximal order,
\item being Auslander-Gorenstein,
\item having finite global dimension,
\item having finite Krull dimension.
\end{enumerate}
\end{thm}
\begin{proof} See \cite{LVO}.
\end{proof}

We immediately obtain from Theorems (\ref{GrR}) and (\ref{lift}), and Corollary (\ref{PBW}), the

\begin{cor} If $G$ is a uniform pro-$p$ group, then $\Omega_G$ is a noetherian, Auslander-regular, scalar local domain which is a maximal order in its quotient division ring of fractions.
\end{cor}

\section{Extensions over finite index}\label{extend} For an arbitrary $p-$adic analytic group $G$, many fundamental
properties of $\Omega_G$ (and of $\Lambda_G$) can be analysed using Corollary (\ref{lift}) and
(\ref{iwa})(\ref{crossed}).

\subsection{Complete noetherian (semi)local rings} \label{local} Recall that a ring $R$ is \emph{semilocal}
if the factor of $R$ by its Jacobson radical $J(R)$ is semisimple artinian. It is \emph{local} if $R/J(R)$ is
simple artinian, and \emph{scalar local} if $R/J(R)$ is a division ring. For a crossed product $R = S \ast H$ of a
finite group $H$, like that in (\ref{iwa})(\ref{crossed}), it's not hard to show that $J(S) \subseteq J(R)$,
\cite[Theorem 1.4.2]{Pas}. From this, Theorem (\ref{comp}) and Corollary (\ref{lift}), and their analogues for
$\Lambda_G$, we deduce (1) of the following. Both it and (2) were known to Lazard.

\begin{thm}
Let $G$ be a compact $p-$adic analytic group.
\begin{enumerate}
\item $\Omega_G$ and $\Lambda_G$ are complete noetherian semilocal
rings. \item $\Omega_G$ and $\Lambda_G$ are (scalar) local rings
if and only if $G$ is a pro$-p$ group.
\end{enumerate}
\end{thm}

\subsection{Primeness and semiprimeness} \label{sprime} The characterisations of
these properties given in the theorem below exactly parallel the results for ordinary group algebras proved in the
early 1960s by I.G. Connell and D.S. Passman \cite[Theorems 4.2.10 and 4.2.14]{Pas2}. However, the proofs here are
quite different from the classical setting; that the stated conditions are necessary is easy to see, but
sufficiency in (1) and (2) depends on Corollary (\ref{lift}) to handle the uniform case, together with non-trivial
results on crossed products of finite groups. Part (3) is much easier - one can simply appeal to the fact (a
consequence of Maschke's theorem) that the group ring of a finite group over a commutative coefficient domain of
characteristic 0 is semiprime, together with the fact that, by definition, $\Lambda_G$ is an inverse limit of such
group rings.

\begin{thm}  Let $G$ be a compact $p-$adic analytic group. \begin{enumerate}
\item \cite{AB} $\Omega_G$ and $\Lambda_G$ are prime if and only
if $G$ has no non-trivial finite normal subgroups. \item \cite{AB} $\Omega_G$ is semiprime if and only if $G$ has
no non-trivial finite normal subgroups of order divisible by $p$. \item $\mathrm{(Neumann,}$ \cite{N}$\mathrm{)}$
$\Lambda_G$ is always semiprime.
\end{enumerate}
\end{thm}

\subsection{Zero divisors} \label{zerodiv} There is a method, familiar from the treatment of ordinary group rings,
which allows one to use homological properties to deduce results about the non-existence of zero divisors in
certain noetherian rings. In its simplest form, which is all that is needed here, the statement is due to Walker
\cite{W}: \emph{if $R$ is a scalar local noetherian semiprime ring of finite global dimension, then $R$ is a
domain.}\footnote{It is a famous and long-standing open question in ring theory whether ``semiprime" is necessary
in Walker's theorem.} This yields the following result; it was proved by Neumann \cite{N} for $\Lambda_G$, but for
$\Omega_G$ it was necessary to wait first for semiprimeness to be settled, as in Theorem (\ref{sprime})(2).

\begin{thm}
Let $G$ be a compact $p-$adic analytic group. Then $\Omega_G$ and $\Lambda_G$ are domains if and only if $G$ is
torsion free.
\end{thm}

\begin{proof}
If $1 \neq x \in G$ with $x^n = 1$, then $(1-x)(1+x+ \cdots x^{n-1}) = 0,$ so the absence of torsion is clearly
necessary. Suppose that $G$ is torsion free. Since $G$ has a pro$-p$ subgroup of finite index by (\ref{uniform}),
its Sylow $q-$subgroups are finite for primes $q$ not equal to $p$. Since $G$ is torsion free these subgroups are
trivial, so $G$ is a pro$-p$ group. Therefore $\Omega_G$ and $\Lambda_G$ are scalar local and noetherian by
Theorem (\ref{local}). The other conditions needed for Walker's theorem are given by Theorems (\ref{sprime})(2) and
(3) and (\ref{global}).
\end{proof}

\subsection{Maximal orders} \label{max} It might seem natural to
suppose, in the light of Theorem (\ref{lift})(3), that whenever
$\Lambda_G$ or $\Omega_G$ are prime then they are maximal orders.
This guess is wrong, though, as the following example shows.

\bigskip

\noindent \textbf{Example:} Let $D := A\rtimes \langle \gamma
\rangle,$ where $A$ is a copy of $\mathbb{Z}_2$ and $\gamma$ is
the automorphism of order 2 sending each 2-adic integer to its
negative. Since $D$ is a pro-2 group with no non-trivial finite
normal subgroups, Theorems (\ref{local}) and (\ref{sprime}) show that
$\Omega_D$ and $\Lambda_D$ are prime noetherian scalar local
rings. But it's not hard to see that neither of these algebras is
a maximal order: for $\Omega_D$, observe that it is local with
reflexive Jacobson radical $J$ which is not principal, impossible
for a prime noetherian maximal order by \cite[Th\'eor\`eme
IV.2.15]{MaR}; for $\Lambda_D,$ the kernel of the canonical map to
$\Zp$ is a reflexive prime ideal which is not localisable,
impossible in a maximal order by \cite[Corollaire IV.2.14]{MaR}. We therefore ask:

\begin{quest} \label{MaxOrd} When are $\Omega_G$ and $\Lambda_G$ maximal orders?
\end{quest}

\noindent Since the powerful structural results \cite{Ch} which can be obtained for certain quotient categories of the
category of finitely generated modules over a noetherian maximal order are potentially important tools in
arithmetic applications \cite{CSS}, this question is of more than passing interest.

In the next three paragraphs we offer a conjecture for the answer to Question \ref{MaxOrd}, and give some evidence in its
support.

\subsection{Conjectured answer to Question \ref{MaxOrd}}

We will need some group-theoretic notions. Let $H$ be a closed subgroup of a compact $p-$adic analytic group $G$.
We say that $H$ is \emph{orbital} if $H$ has finitely many $G-$conjugates, or equivalently if the normaliser $N$
has finite index in $G$. We say that an orbital subgroup $H$ is \emph{isolated} if $N/H$ has no non-trivial finite
normal subgroups.

We will say that $G$ is \emph{dihedral-free} if, whenever $H$ is an orbital closed subgroup of $G$ with $\dim H =
1$, $H \cong \Zp$. This seems to be the correct generalisation of the definition in \cite{Br}.

\begin{conj} Let $G$ be a compact $p-$adic analytic group, and suppose $\Omega_G$ is prime. Then
$\Omega_G$ is a maximal order if and only if $G$ is dihedral-free.
\end{conj}

\subsection{Necessary conditions on $G$}

We fix a prime $p$ and assume throughout this paragraph that $G$ is a compact $p-$adic analytic group.

\begin{prop} Suppose $\Omega_G$ is a prime maximal order and let $H$ be a closed \emph{normal} subgroup of $G$ with $\dim H = 1.$
 Then $H$ is pro-$p$.
\end{prop}
\begin{proof}
We may assume that $H$ is isolated, so $G/H$ has no non-trivial finite normal subgroups. Hence, by Theorem
(\ref{sprime})(1), $w_H = \ker(\Omega_G \to \Omega_{G/H})$ is a prime ideal of $\Omega_G$, and it is not hard to see
that it is also a reflexive ideal.\footnote{One quick way to see this uses the canonical dimension from
(\ref{Cdim}): since $\Cdim(\Omega_G / w_H) = \dim (G/H) = \dim G - 1$ and since $\Omega_G$ is
Auslander-Gorenstein, $w_H$ is reflexive by Gabber's Maximality Principle \cite[Theorem 2.2]{Staf}.} Now because
$\Omega_G$ is a maximal order and $w_H$ is a prime reflexive ideal, it must be localisable \cite[Corollaire
IV.2.14]{MaR}.

But the conditions needed for augmentation ideals to be localisable are known \cite[Theorem E]{AB}: $H/F$ must be
pro-$p,$ where $F$ is the largest finite normal $p'-$subgroup of $H$. Since $H$ is normal in $G$ and $G$ has no
non-trivial finite normal subgroups, $F = 1$ and $H$ is pro-$p$ as required.
\end{proof}

We need the following  group-theoretic lemma. We first set
$\epsilon$ to be 1 for $p$ odd, and $\epsilon = 2$ if $p=2$, and
define, for a closed normal uniform subgroup $N$ of $G,$ $E_G(N)$
to be the centraliser in $G$ of $N/N^{p^\epsilon},$
\cite[(2.2)]{AB}.

\begin{lem} Suppose that $G$ is a pro-$p$ group of finite rank with no non-trivial
finite normal subgroups. Let $N$ be a maximal open normal uniform subgroup of $G$. Then
\[E_G(N) = N.\]
\end{lem}
\begin{proof}Recall that $E = E_G(N)$ is an open normal subgroup of $G$ containing $N$.
If $E$ strictly contains $N$ then $E/N$ must meet the centre $Z(G/N)$ non-trivially since $G/N$ is a finite
$p-$group by \cite[Proposition 1.11(ii)]{DDMS}. Pick $x \in E \backslash N$ such that $xN \in Z(G/N)$; then $H =
\langle N,x\rangle$ is normal in $G$ by the choice of $x$, and also $H$ is uniform by \cite[Lemma 2.3]{AB}. This
contradicts the maximality of $N$.
\end{proof}

\begin{cor} Let $H$ be a pro-$p$ group of finite rank with no non-trivial finite
normal subgroups. Suppose that $\dim H = 1$. Then $H \cong \Zp$, unless $p=2$ and $H$ is isomorphic to $D$.
\end{cor}
\begin{proof} Choose a maximal open normal uniform subgroup $N$ of $H$.
By the proposition and the lemma, $H/N \hookrightarrow \Aut(N/N^{p^\epsilon})$. If $p$ is odd, $|N :
N^{p^\epsilon}| = p,$ so the latter automorphism group is just $\Fp^\times$. Since $H/N$ is a $p-$group by
\cite[Proposition 1.11(ii)]{DDMS} again, $H = N \cong \Zp$. If $p=2$ and $H > N$, $H \cong D$.
\end{proof}

This gives us the following weak version of one half of the conjecture. To improve the result from ``normal'' to ``orbital'' will presumably require some technical work on induced ideals.

\begin{cor} Suppose $\Omega_G$ is a prime maximal order. Then any closed normal subgroup $H$ of $G$ of dimension 1
is isomorphic to $\Zp$.
\end{cor}

\begin{proof}
When $p$ is odd the statement is immediate from the proposition and corollary above. So suppose that $p = 2.$ We have to rule out the possibility that $H \cong D$, so suppose for a contradiction that this is the case. Then, as in the proof of the proposition, $w_H$ is a prime reflexive, and hence localisable, ideal of $\Omega_G.$ Let $R$ denote the local ring $(\Omega_G)_{w_H},$ which has global dimension one by \cite[Th\'eor\`eme IV.2.15]{MaR}. Let $C = \langle c \rangle$ be a copy of the cyclic group of order 2 in $H$. Then $\mathbb{F}_2 C \subseteq \Omega_G$ and $\Omega_G$ is a projective $\mathbb{F}_2 C-$module by \cite[Lemma 4.5]{Bru}. Thus $R$ is a flat $\mathbb{F}_2 C-$module. Since $c + 1 \in J(R),$ the $\mathbb{F}_2
C-$module $R/J(R)$ is a sum of copies of the trivial module, so
\[ \infty = \pd_{\mathbb{F}_2 C} (\mathbb{F}_2) = \pd_{\mathbb{F}_2 C} (R/J(R)) \leq \pd_R (R/J(R)) = 1. \]
This contradiction shows that the only possibility for $H$ is $\mathbb{Z}_2.$
\end{proof}

\subsection{Sufficient conditions on $G$}
\label{Suff} We use the following result, essentially due to R. Martin:

\begin{prop}\cite{Martin} Let $R$ be a prime noetherian maximal order and let $F$ be a finite group. Let $S = R \ast F$ be a
prime crossed product. Then $S$ is a maximal order if and only if
\begin{enumerate}[{(}a{)}]
\item every reflexive height 1 prime $P$ of $S$ is localisable, and
\item $\gld(S_P) < \infty$ for all such $P$.
\end{enumerate}
\end{prop}
\begin{proof} Conditions (a) and (b) hold in any prime noetherian maximal order, \cite[Th\'eor\`eme IV.2.15]{MaR}.
Conversely, suppose that (a) and (b) hold. We use the Test Theorem \cite[Theorem 3.2]{Martin}. Condition (i) of
the Test Theorem is just condition (a). We claim that if $P$ is as in the theorem, then $\gld(S_P) = 1.$ It's
easy to check that $P \cap R$ is a semiprime reflexive ideal of $R$, so that the localisation $R_{P\cap R}$ exists
and is hereditary by \cite[Th\'eor\`eme IV.2.15]{MaR}. Thus $R_{P\cap R} \ast F$ has injective
dimension 1 by (\ref{ausgor})(1). But $S_P$ is a localisation of $R_{P\cap R} \ast F,$ so - given (b) and the comments in (\ref{lift}) - $\gld
(S_P) \leq 1.$ The reverse inequality is obvious, so our claim follows. Condition (ii) now follows from \cite[Proposition 2.7]{Martin}. Condition (iii) follows from the proof of \cite[Lemma 3.5]{Martin} and condition (iv) follows from \cite[Remark 3.6 and Lemma 3.7]{Martin}.
\end{proof}

\begin{lem} Let $G$ be a pro-$p$ group of finite rank with no non-trivial finite normal subgroups.
Then every reflexive height 1 prime of $\Omega_G$ is localisable.
\end{lem}
\begin{proof}
Let $P$ be a reflexive height 1 prime of $\Omega_G$. Choose an open normal uniform subgroup $N$ of $G.$ Then
$\Omega_N$ is a maximal order by Corollary (\ref{lift}). Set $\overline{G} := G/N.$ Now let $Q = P \cap \Omega_N$ -
it is easy to see \cite[Remark 3.6]{Martin} that this is a height 1 reflexive $\overline{G}-$prime ideal of
$\Omega_N$. Indeed, $Q$ is the intersection of a $\overline{G}-$orbit of reflexive prime ideals $\{P_1, \ldots ,
P_n \}$ of $\Omega_N.$

Since each $P_i$ is localisable by \cite[Th\'eor\`eme IV.2.15]{MaR}, $Q$ is localisable. In other
words, the subset $\mathcal{C} := \mathcal{C}_{\Omega_N} (Q) = \cap_{i=1}^n \mathcal{C}_{\Omega_N}(P_i)$ is a
$\overline{G}-$invariant Ore set in $\Omega_N$. An easy calculation \cite[proof of Lemma 13.3.5(ii)]{Pas2} shows
that $\mathcal{C}$ is an Ore set in $\Omega_G$. In other words, the semiprime ideal $A = \sqrt{Q\Omega_G}$ is
localisable in $\Omega_G$ and
\[(\Omega_N)_Q \ast \overline{G} \cong (\Omega_G)_A.\]
Since $\overline{G}$ is a $p-$group, $A = P$ by \cite[Proposition 16.4]{Pas} and the result follows.
\end{proof}

\begin{cor} Let $G$ be a torsion free compact $p-$adic analytic group. Then $\Omega_G$ is a prime maximal order.
\end{cor}
\begin{proof} Suppose that $G$ is as stated. Since $G$ has a pro$-p$ open subgroup,
 the Sylow $q-$subgroups of $G$ are finite, and hence trivial, for all primes $q$ not equal to $p$. That is, $G$ is a pro$-p$ group.
 Thus the corollary  follows from the lemma and the proposition, since $\gld \Omega_G$ is finite by Theorem (\ref{global}).
\end{proof}

\section{Dimensions}
\label{Dims}

\subsection{Global dimension}\label{global} The situation as regards the global dimension of $\Omega_G$ and
$\Lambda_G$ is completely understood, and depends fundamentally on properties of the cohomology of profinite
groups - in particular behaviour under finite extensions - due to Serre \cite{Serre}. The result is due to Brumer
\cite[Theorem 4.1]{Bru} who computed the global dimension of the completed group algebra of an arbitrary profinite
group $G$ with coefficients in a pseudo-compact ring $R$. As a consequence of his work, we have

\begin{thm} Let $G$ be a compact $p-$adic analytic group of dimension $d.$ Then $\Omega_G$ and $\Lambda_G$ have
 finite global dimension if and only if $G$ has no elements of order $p,$ and in this case
 $$ \gld (\Omega_G) = d \quad \textit{ and } \quad \gld (\Lambda_G) = d+1. $$
\end{thm}

\subsection{Auslander-Gorenstein rings} \label{ausgor} Recall that the group algebra of an arbitrary finite group
over any field is a Frobenius algebra \cite[Proposition 4.2.6]{Wei}, and thus is self-injective. It should therefore come as no
surprise that injective dimension is well-behaved for Iwasawa algebras. In fact, much more is true:

\begin{thm}
\cite[Theorem J]{AB} Let $G$ be a compact $p-$adic analytic group of dimension $d.$ Then $\Omega_G$ and
$\Lambda_G$ are Auslander-Gorenstein rings of dimensions $d$ and $d+1$ respectively.
\end{thm}

This result was first proved by O. Venjakob \cite{Ven1} and is easy to deduce from Theorems (\ref{lift})(4) and (\ref{global}), as follows. Let $H$ be an open uniform normal subgroup of $G.$ Then $\Omega_H$ and $\Lambda_H$ are Auslander-Gorenstein by Theorem (\ref{lift})(4), and the dimensions are given by Theorem (\ref{global}). Now apply (\ref{iwa})(\ref{crossed}): a simple lemma \cite[Lemma 5.4]{AB} shows that
\begin{equation} \Ext^i_{\Omega_G}(M, \Omega_G) \cong \Ext^i_{\Omega_H}(M,\Omega_H) \tag{1}\label{frob}  \end{equation}
for all $i \geq 0$ and all $\Omega_G -$modules $M$, with a similar isomorphism for $\Lambda_G$, and the result
follows.

\subsection{Dimension functions for Auslander-Gorenstein rings}\label{dim} We recall from \cite{Le} the basics of dimension theory
over an Auslander-Gorenstein ring $R$. Write $d$ for the injective dimension of $R$. The \emph{grade} $j(M)$ of a
finitely generated $R-$module $M$ is defined as follows:
\[j(M) = \min\{j : \Ext^j_R(M, R) \neq 0\}.\]
Thus $j(M)$ exists and belongs to the set $\{0,\ldots , d \} \cup \{+\infty \}$. The \emph{canonical dimension} of
$M$, $\Cdim(M)$ is defined to be
\[\Cdim(M) = d - j(M).\]
It is known \cite[Proposition 4.5]{Le} that $\Cdim$ is an exact, finitely partitive dimension function on finitely
generated $R-$modules in the sense of \cite[\S 6.8.4]{MR}. That is,
\begin{itemize}
\item $\Cdim(0) = - \infty;$ \item if $0 \longrightarrow N
\longrightarrow M \longrightarrow T \longrightarrow 0$ is an exact
sequence of finitely generated modules, then $\Cdim (M) =
\textrm{max} \{\Cdim (N), \Cdim (T) \}$; \item if $MP = 0$ for a
prime ideal $P$ of $R$, and $M$ is a torsion $R/P-$module, then
$\Cdim (M) \leq \Cdim (R/P) - 1;$ \item if $\Cdim (M) = t$ then
there is an integer $n$ such that every descending chain $M = M_0
\supseteq M_1 \supseteq \cdots \supseteq M_i \supseteq M_{i+1}
\cdots $ of submodules of $M$ has at most $n$ factors
$M_i/M_{i+1}$ with $\Cdim (M_i/M_{i+1}) = t.$
\end{itemize}
The ring $R$ is said to be \emph{grade symmetric} if
\[\Cdim(_RM) = \Cdim(M_R)\]
for any $R$-$R-$bimodule $M$ which is finitely generated on both sides.\footnote{Alternatively, we can say in these
circumstances that the dimension function $\Cdim$ is \emph{symmetric}.} The triangular matrix ring $\begin{pmatrix}k & k \\
0 & k\end{pmatrix}$ over a field $k$ gives an easy example of an Auslander Gorenstein ring which is \emph{not}
grade symmetric.

The existence of an exact, finitely partitive, symmetric dimension function for the finitely generated modules
over a noncommutative noetherian ring $R$ is a very valuable tool which is often not available: the
Gel'fand-Kirillov dimension \cite[\S 8.1]{MR} - although symmetric - is often not defined; and although the Krull
dimension is always defined \cite[\S 6.2]{MR}, it is a long-standing open question whether it is symmetric in
general. As we shall see in the next paragraph, the canonical dimension function fulfils these requirements for an
Iwasawa algebra.

If $\delta$ is a dimension function on finitely generated
$R-$modules, we say that $R$ is \emph{Cohen-Macaulay with respect
to $\delta$} if $\delta(M) = \Cdim(M)$ for all finitely generated
$R-$modules $M$.

This definition is consistent with, and therefore generalises, the
definition from commutative algebra. To see this, suppose that $R$
is a commutative noetherian ring of locally constant dimension
$d.$ Suppose that $R$ is Cohen-Macaulay \cite[Definition
2.1.1]{BH}, and let $M$ be a finitely generated $R-$module with
Krull dimension $\mathcal{K}(M)$, (which equals the
Gel'fand-Kirillov dimension of $M$ if $R$ is affine). Then
\begin{equation} j(M) + \mathcal{K}(M) = d, \tag{1}\label{comm}
\end{equation}
\cite[Corollary 2.1.4 and Theorem 1.2.10(e)]{BH}. And conversely,
if (\ref{comm}) holds for all simple $R-$modules $M$, then $R$ is
Cohen-Macaulay \cite[Theorem 1.2.5]{BH}.

\subsection{Canonical dimension for $\Omega_G$}
\label{Cdim} We continue in this paragraph to assume that $G$ is a
compact $p-$adic analytic group of dimension $d$. Fix an open
uniform normal subgroup $H$ of $G,$ and let $M$ be a finitely
generated $\Omega_G-$module. By Theorem (\ref{ausgor}) and paragraph
(\ref{dim}), and with the obvious notation, $\Cdim_G(-)$ and
$\Cdim_H(-)$ are well-defined dimension functions, and in fact
(\ref{ausgor})(\ref{frob}) shows that \begin{equation}
\label{index} \Cdim_H(M) = \Cdim_G(M).  \tag{1}\end{equation} In
particular, in studying the canonical dimension we may as well
assume that $G = H$ is uniform, which we now do. Hence, by Theorem
(\ref{GrR}), the graded ring of $\Omega_G$ is a polynomial
$\Fp-$algebra in $d$ variables.

Choose a good filtration for $M$ ($F_nM = MJ^{-n}$ for $n \leq 0$
will do) and form the associated graded module $\gr M$. Because
the $J-$adic filtration is Zariskian, it follows from \cite[Remark
5.8]{BjE} that \begin{equation} j(\gr M) = j(M).\tag{2}
\label{fungr} \end{equation} Moreover, from this and the
concluding remarks  of (\ref{dim}) we see that
\begin{equation}
\label{krull} \Kr(\gr M) = \Cdim(\gr M) = d - j(M).\tag{3}
\end{equation}

(This shows, incidentally, that $\Kr(\gr M)$ is actually independent of the choice of
good filtration on $M$.)\footnote{Consider (\ref{krull}) with $M$ the trivial $\Omega_G-$module $\Fp$.
Then $\Kr(\gr M) = 0,$ so $j(M) = d$ and therefore the injective dimension of $\Omega_G$ actually
equals $d$, providing another proof of the numerical part of Theorem (\ref{global}).} Combining (\ref{fungr}) and (\ref{krull}), we find that
\[\Cdim(M) = d - j(M) = \Cdim(\gr M) = \Kr(\gr M) = \GK(\gr M)\]
for any choice of good filtration on $M$. This proves the last
part of the

\begin{prop} Let $G$ be a compact $p-$adic analytic group.
\begin{enumerate}
\item $\Omega_G$ is grade-symmetric.
\item $\Omega_G$ is ideal-invariant with respect to $\Cdim$.
\item Suppose that $G$ is uniform. Then for all finitely generated $\Omega_G-$modules $M$,
\[\Cdim(M)\quad  = \quad  \GK(\gr M).\]
\end{enumerate}
\end{prop}
\begin{proof}
(1) In view of (\ref{Cdim})(\ref{index}) we can and do assume that
$G$ is uniform. Write $J$ for the Jacobson radical of $\Omega_G$
and let $M$ be a finitely generated $\Omega_G-$module. Then by the
definition of the Gel'fand Kirillov dimension \cite[\S
8.1.11]{MR}, $\GK(\gr M)$ is the growth rate $\gamma(f)$ of the
function
\[f(n) = \dim \frac{M}{MJ^n};\]
note that this function is eventually polynomial because the finitely
generated $\gr \Omega_G-$module $\gr M$ has a Hilbert polynomial.

Now let $N$ be an $\Omega_G-$bimodule, finitely generated on both sides.
Then $NJ$ is a sub-bimodule, and $N/NJ$ is finite dimensional over $\Fp$ because $N$ is a finitely
generated right $\Omega_G-$module. Hence $N/NJ$ is also a finite dimensional \emph{left} $\Omega_G-$module and
as such is killed by some power of $J$, $J^a$ say. Thus $J^aN \subseteq NJ$ and similarly there exists an
integer $b\geq 1$ such that $NJ^b \subseteq JN.$ An easy induction on $n$ shows that
\begin{equation}
\label{Compare} J^{abn}N \subseteq NJ^{bn} \subseteq J^nN
\end{equation}
for all $n\geq 0$.
Letting $f(n) = \dim \frac{N}{NJ^n}$ and $g(n) = \dim \frac{N}{J^nN}$, we obtain
\[g(n) \leq f(bn) \leq g(abn)\]
for all $n\geq 0$. It follows that $\Cdim(N_{\mid \Omega_G}) =
\gamma(f) = \gamma(g) = \Cdim(_{\Omega_G \mid}N)$, proving part
(1).

For part (2), recall \cite[6.8.13]{MR} that a ring $R$ is said to
be \emph{ideal-invariant} with respect to a dimension function
$\delta$ if $\delta(M \otimes_R I) \leq \delta(M)$ for all
finitely generated right $R-$modules $M$ and all two-sided ideals
$I$ of $R$ and if the left-hand version of this statement also
holds.

In fact, we will show that \begin{equation}\label{idinv} \Cdim(M
\otimes_{\Omega_G} N) \leq \Cdim(M)\tag{4}\end{equation} for any
finitely generated $\Omega_G-$module $M$ and any
$\Omega_G-$bimodule $N$, finitely generated on both
sides.\footnote{Compare this with \cite[Proposition 8.3.14]{MR}.}
Let $M$ and $N$ be as above, and let $H$ be an open uniform normal
subgroup of $G$. Since there is an $\Omega_H-$epimorphism $M
\otimes_{\Omega_H} N \twoheadrightarrow M \otimes_{\Omega_G} N$,
(\ref{ausgor})(\ref{index}) shows that we can replace $G$ by $H$
in proving (\ref{idinv}); that is, we now assume that $G$ is
uniform.

Choose the integer $a$ as above so that $J^{an}N \subseteq NJ^n$ for all $n\geq 0$. Fix $n$ and let
\[f(n) = \dim \frac{M}{MJ^n}\quad\mbox{and}\quad g(n) = \dim \left(\frac{M\otimes_{\Omega_G} N}{(M\otimes_{\Omega_G}
 N).J^n}\right).\]
Note that $(M\otimes_{\Omega_G} N).J^n$ equals the image of $M
\otimes_{\Omega_G} NJ^n$ in $M\otimes_{\Omega_G} N$ so the
right-exactness of tensor product gives
\[M \otimes_{\Omega_G} \left(\frac{N}{J^{an} N}\right) \twoheadrightarrow M \otimes_{\Omega_G}
 \left(\frac{N}{NJ^n}\right) \cong \frac{M\otimes_{\Omega_G} N}{(M\otimes_{\Omega_G} N).J^n}.\]
Now we have a natural isomorphism of right $\Omega_G-$modules
\[M \otimes_{\Omega_G} \frac{N}{J^{an} N} \cong \frac{M}{MJ^{an}} \otimes_{\Omega_G} N\]
and picking a finite generating set of size $t$ for the left
$\Omega_G-$module $N$ shows that
\[\dim \left(\frac{M}{MJ^{an}} \otimes_{\Omega_G} N\right) \leq \left(\dim \frac{M}{MJ^{an}}\right) \cdot t.\]
Hence
\[g(n) = \dim \left(\frac{M\otimes_{\Omega_G} N}{(M\otimes_{\Omega_G} N).J^n}\right) \leq \dim \left(M \otimes_{\Omega_G} \left(\frac{N}{J^{an} N}\right)\right) \leq f(an)\cdot t\]
for all $n\geq 0$, so $\Cdim(M\otimes_{\Omega_G} N) = \gamma(g)
\leq \gamma(f) = \Cdim(M)$ as required.
\end{proof}

The above proposition is due to the first author; it was inspired by a result of S. J. Wadsley \cite[Lemma
3.1]{Wadsley}.

\subsection{Characteristic varieties} \label{char}
Assume in this paragraph that $G$ is uniform. Let $M$ be a
finitely generated $\Omega_G -$module. There is another way of
seeing that $\Kr(\gr M)$ does not depend on the choice of good
filtration for $M$, as follows. It is well known \cite[Chapter III, Lemma 4.1.9]{LVO}
that
\[J(M) := \sqrt{\Ann_{\gr \Omega_G}\left(\gr M\right)}\]
is independent of this choice. Standard commutative algebra now gives
\[\Kr(\gr M) = \Kr\left(\frac{\gr \Omega_G}{J(M)}\right),\]
as claimed.

The graded ideal $J(M)$ is called the \emph{characteristic ideal} of $M$, and the affine variety $\Ch(M)$ defined by it is called the \emph{characteristic variety} of $M$. Thus we obtain yet another expression for the canonical dimension of $M$:
\begin{equation}
\label{ChDim}
\Cdim(M) = \dim \Ch(M).
\end{equation}

The characteristic variety is defined in an entirely analogous fashion for finitely generated modules over enveloping algebras and Weyl algebras. In that setting it enjoys many pleasant properties, in addition to the simple formula (\ref{ChDim}). In particular, there exists a \emph{Poisson structure} on $\Ch(M)$, which gives more information about $M$ through the geometric properties of the characteristic variety. For example, the fact that the characteristic variety of a finitely generated $A_n(\mathbb{C})-$module is integrable can be used to prove the Bernstein inequality.

\begin{quest} Is there a way of capturing more information about $M$ in the characteristic variety $\Ch(M)$?
\end{quest}

The naive method (mimicking the construction of the Poisson structure in the enveloping algebra case) seems to fail because derivations are not sufficient when studying algebras in positive characteristic: they kill too much. Presumably, if the answer to the above question is affirmative, then differential operators in characteristic $p$ will play a role.

\subsection{No GK-dimension} The theory outlined in the previous sections will sound very familiar to the experts. However, Iwasawa algebras are \emph{not} Cohen Macaulay with respect to the GK dimension. This is easily seen by decoding the definition of GK dimension in the case when $G \cong \Zp$: in this case, $\Omega_G$ is isomorphic to the one-dimensional power series ring $\Fp[[t]]$, which (being uncountable) contains polynomial algebras over $\Fp$ of arbitrarily large dimension. Thus $GK(\Omega_G) = \infty$ for any infinite $G$, since any such $G$ will contain a closed subgroup isomorphic to $\Zp$.

If one tries to brush this problem away by replacing the GK dimension by the canonical dimension, then one has to be careful not to fall into the following trap.

Recall \cite[Lemma 8.1.13(ii)]{MR} that if $R \subseteq S$ are affine
$k-$algebras over a field $k$, then for any finitely generated
$S-$module $M$,
\begin{equation}
\label{niceGK}
\GK(N) \leq \GK(M)
\end{equation}
whenever $N$ is a finitely generated $R-$submodule of $M$. This
enables one to ``pass to subalgebras of smaller dimension" and use
inductive arguments on the GK dimension - a ploy used, for
example, in the computation of the Krull dimension of
$\mathcal{U}(\mathfrak{sl}_2(\mathbb{C}))$ by S.P. Smith \cite[Theorem 8.5.16]{MR}. Another consequence of this property of GK dimension is
that it is impossible to find an embedding $R \hookrightarrow S$
of $k-$algebras such that $\GK(R) > \GK(S)$.

Unfortunately, (\ref{niceGK}) fails for Iwasawa algebras, if one tries to replace the GK dimension by the canonical dimension. This is due to the following pathological example:

\begin{example} \cite[Chapter VII, page 219]{SZ} There exists a continuous embedding of $\Fp-$algebras
\[\Omega_G \hookrightarrow \Omega_H\]
where $\dim G = 3$ and $\dim H$ = 2.
\end{example}
\begin{proof}
Let $G = \Zp^3$ and $H = \Zp^2$. By Theorem (\ref{PBW}) we can identify $\Omega_G$ with the three-dimensional power series ring $\Fp[[x,y,z]]$ and $\Omega_H$ with the two-dimensional power series ring $\Fp[[a,b]]$.

Because $\Fp[[a]]$ is uncountable, we can find an element $u = u(a) \in a\Fp[[a]]$ such that the $\Fp-$algebra generated by $a$ and $u$ is isomorphic to the two-dimensional polynomial ring $\Fp[a,u]$. Define $\theta : \Fp[[x,y,z]] \to \Fp[[a,b]]$ to be the unique continuous $\Fp-$algebra map such that
\[
\theta(x) = b, \quad \theta(y) = ab, \quad \theta(z) = ub.
\]
We have
\[\theta\left(\sum_{\lambda,\mu,\nu\in\mathbb{N}}r_{\lambda,\mu,\nu}x^\lambda y^\mu z^\nu \right) = \sum_{n=0}^\infty b^n \left(\sum_{\lambda + \mu + \nu = n}r_{\lambda,\mu,\nu} a^\mu u^\nu\right).\]
This shows that $\theta$ is an injection, as required.
\end{proof}

One can of course concatenate these embeddings and produce a continuous embedding of $\Omega_G$ into $\Fp[[a,b]]$ for abelian uniform pro-$p$ groups $G$ of arbitrarily large dimension. Here is the actual counterexample to the analogue of (\ref{niceGK}).

\begin{example} There exist uniform pro-$p$ groups $H \subset G$, a finitely generated $\Omega_G-$module
$M$ and a finitely generated $\Omega_H-$submodule $N$ of $M$ such that $\Cdim(M) = 2$, but $\Cdim(N) = 3$.
\end{example}
\begin{proof}Let $R = \Fp[[a,b,c,d]]$ and $S = \Fp[[b,c,d]]$. Let $I$ be the ideal of $R$ generated by $c - ab$ and $d - u(a)b$ where $u(a)$ is chosen as in the previous example and let $M = R/I$. By construction, the graded ideal $\gr I$ is generated by the symbols of $c$ and $d$, so
\[\Cdim(M) = \Kr(\gr M) = 2.\]
Now if $r \in I \cap S$, then $\theta(r) = 0$, letting $\theta: \Fp[[b,c,d]] \hookrightarrow \Fp[[a,b]]$ be as above. Hence $r = 0$, so $S \hookrightarrow R/I = M$. Therefore the cyclic $S-$submodule $N$ of $M$ generated by $1 + I$ is actually free, so $\Cdim(N) = 3$.
\end{proof}

\subsection{Krull dimension}
\label{KdimOmega}
The Krull-(Gabriel-Rentschler) dimension of $\Omega_G$ was first studied by one of the authors in \cite{Ard}. An immediate upper bound of $\dim G$ can be obtained using Theorem (\ref{lift}), or if one prefers, using \cite[Corollary 1.3]{ASZ}. Here is a result covering a large number of cases.

\begin{thm} \cite[Theorem A and Corollary C]{Ard}
Let $G$ be a compact $p-$adic analytic group, and let
$\mathfrak{g}$ be the $\Qp -$Lie algebra of an open uniform
subgroup of $G.$ Let $\mathfrak{r}$ denote the soluble radical of
$\mathfrak{g}$ and suppose that the semisimple part
$\mathfrak{g}/\mathfrak{r}$ of $\mathfrak{g}$ is a direct sum of
some number of copies of $\mathfrak{sl}_2(\Qp)$. Then
\[\Kr(\Omega_G) = \dim G.\]
\end{thm}

In particular, $\Kr(\Omega_G)$ equals $\dim G$ whenever $G$ is
soluble-by-finite, and whenever $\dim G \leq 3$. The main idea in
the proof is to obtain a lower bound on the Krull dimension of
$\Omega_G$ for \emph{any} compact $p-$adic analytic group $G$.
Namely, with $\mathfrak{g}$ as in the theorem, and writing
$\lambda(\mathfrak{g})$ for the length of the longest chain of
subalgebras of $\mathfrak{g}$, we have $$\lambda(\mathfrak{g}) \leq
\Kr(\Omega_G).$$

\begin{quest}
\label{Kdim}
With the above notation, is $\Kr(\Omega_G) = \lambda(\mathfrak{g})$ in general?
\end{quest}

It is easy to see that $\lambda(\mathfrak{g}) = \lambda(\mathfrak{n}) + \lambda(\mathfrak{g}/\mathfrak{n})$ whenever $\mathfrak{n}$ is an ideal of $\mathfrak{g}$. Let $N$ be a closed uniform subgroup of $G$ with Lie algebra $\mathfrak{n}$.

\begin{quest}
\label{Breakup}
Is $\Kr(\Omega_G) = \Kr(\Omega_N) + \Kr(\Omega_{G/N})$?
\end{quest}

Aside from its intrinsic interest, an affirmative answer to
Question \ref{Breakup} would obviously reduce Question \ref{Kdim} to the study of
\emph{almost simple} groups $G$, (where we say that a uniform
pro-$p$ group $G$ is \emph{almost simple} provided its Lie algebra has no
non-trivial ideals).

The classical split simple Lie algebras are the first examples to study. Given such a Lie algebra $\mathfrak{g}$, choose a Borel subalgebra $\mathfrak{b}$ and a Cartan subalgebra $\mathfrak{t}$. Then it is easy to produce a chain of subalgebras of $\mathfrak{g}$ of length $\dim \mathfrak{b} + \dim \mathfrak{t}$.

\begin{quest}
\label{KdimAlmostSimple}
For $G$ almost simple and split, is $\Kr(\Omega_G) = \dim \mathfrak{b} + \dim \mathfrak{t}$ ?
\end{quest}

Question \ref{KdimAlmostSimple} has an affirmative answer in the two smallest cases:
$\mathfrak{g} = \mathfrak{sl}_2(\Qp)$ and $\mathfrak{g} =
\mathfrak{sl}_3(\Qp)$. In particular,

\begin{thm}\cite[Theorem B]{Ard}.
Let $G$ be a uniform pro-$p$ group with $\Qp-$Lie algebra
$\mathfrak{sl}_3(\Qp)$. Then $\Omega_G$ is a scalar local complete
noetherian domain of global dimension 8, with $$\Kr(\Omega_G) =
7.$$
\end{thm}

The main idea of the proof of this last result is to show that $\Omega_G$ has no finitely generated modules whose
canonical dimension equals precisely 1; that is, there is a ``gap" at $\Cdim = 1$.\footnote{A similar idea was
used by Smith \cite{Smith} in giving an upper bound for the Krull dimension of $\mathcal{U}(\mathfrak{g})$ when
$\mathfrak{g}$ is a complex semisimple Lie algebra. We note in passing that $\Kr(\mathcal{U}(\mathfrak{g}))$ when
$\mathfrak{g}$ is complex semisimple has been recently proved to be equal to $\dim \mathfrak{b}$ by Levasseur
\cite{Le2}, answering a long-standing question in the affirmative.} The extra $\dim \mathfrak{t}$ term in our
conjectured formula for $\Kr(\Omega_G)$ comes from the fact that $\Omega_G$ is scalar local - this fact is used
crucially in the proof of the lower bound for the Krull dimension of $\Omega_G$.

\section{Two-sided ideal structure}

\subsection{} One of the first questions asked when studying a noetherian algebra $R$ is ``what are its two-sided
ideals?'' It is usually sensible to focus first on the \emph{prime} ideals of $R$.

One way of answering the above question is to give a reduction to
the commutative case. This is a recurring theme in noncommutative
algebra. For example, if $R = k[G]$ is the group algebra of a
polycyclic group $G$ over a field $k$, the paper \cite{Roseblade}
by J. E. Roseblade achieves this, ``to within a finite
group''.\footnote{See \cite[Chapter 5]{Pas} for more details.}
Similar results hold for universal enveloping algebras
$\mathcal{U}(\mathfrak{g})$ of finite dimensional soluble Lie
algebras over a field $k$: see \cite{Dixm} and \cite[Chapter
13]{MR}. As for the case when $\mathfrak{g}$ is semisimple, one
can view the huge body of research on the primitive ideals of
$\mathcal{U}(\mathfrak{g})$ as an analysis of the failure of the naive hope
that these primitive ideals should be generated by their
intersection with the centre of $\mathcal{U}(\mathfrak{g})$, \cite{Dixm}.
And for quantised function algebras of semisimple groups, and many
related quantum algebras, there are ``stratification theorems''
which describe their prime and primitive spectra as finite
disjoint unions of affine commutative pieces, \cite[Theorem
II.2.13]{BG}.

Unfortunately, no such results are currently known for Iwasawa algebras - see below for a summary of what little
\emph{is} known. Alleviation of this state of gross ignorance would seem to be the most pressing problem in the
subject.

Because of the crossed product decomposition (\ref{iwa})(\ref{crossed}) and the going up and down theorems for
crossed products of finite groups \cite[Theorem 16.6]{Pas}, one should naturally first concentrate on the case
when $G$ is uniform.

\subsection{Ideals arising from subgroups and from centres} \label{ideals}Since centrally generated one-sided ideals are necessarily two-sided, it helps to know the centre of the ring in question. However the centre of Iwasawa algebras is not very big:

\begin{thm}\cite[Corollary A]{Ard2} Let $G$ be a uniform pro$-p$ group and let $Z$ be its centre. Then the centre of $\Omega_G$ equals $\Omega_Z$ and the centre of $\Lambda_G$ equals $\Lambda_Z$.
\end{thm}

Thus when the centre of $G$ is trivial (and this happens frequently), $\Omega_G$ has no non-trivial centrally
generated ideals. This is one place where the analogy with enveloping algebras of semisimple Lie algebras breaks
down.

One can also produce two-sided ideals by using normal subgroups. Certainly when $H$ is a closed normal subgroup of $G$, the augmentation ideal
\[w_H := \ker(\Omega_G \to \Omega_{G/H})\]
is a two-sided ideal of $\Omega_G$ and we can tell whether it is
prime or semiprime using Theorem (\ref{sprime}). As for
$\Lambda_G$, $H$ yields two augmentation ideals: the inverse image
$v_H$ of $w_H$ under the natural projection $\Lambda_G
\twoheadrightarrow \Omega_G$ and ``the'' augmentation ideal
\[I_H = \ker(\Lambda_G \to \Lambda_{G/H}).\]
The behaviour of these ideals regarding localisation is quite well understood:

\begin{thm} Let $H$ be a closed normal subgroup of the compact $p-$adic analytic group $G$ and let $F$ be the largest finite normal subgroup of $H$ of order coprime to $p$. Then
\begin{enumerate}
\item \cite{AB} $w_H$ and $v_H$ are localisable if and only if $H/F$ is pro-$p$,
\item \cite{Ard4} $I_H$ is localisable if and only if $H$ is finite-by-nilpotent.
\end{enumerate}
\end{thm}

Notwithstanding the above, the most embarrassing aspect of the state of our knowledge about ideals of Iwasawa
algebras is the lack of examples. In particular, we've noted that central elements and closed normal subgroups
give rise to ideals. This suggests the following improperly-posed question, for which we'll suggest more precise
special cases in the succeeding paragraphs.

\begin{quest}
Is there a mechanism for constructing ideals of Iwasawa algebras which involves neither central elements nor
closed normal subgroups?
\end{quest}

One way to begin the study of prime ideals is to look first at the
smallest non-zero ones - that is, the prime ideals of height one.
With one eye on the commutative case and another on the results of
(\ref{max}) on maximal orders, one can ask when they are all
principal. Here are two slightly more precise ways to ask this
question:

\begin{quest}
When is $\Omega_G$ a unique factorisation ring in the sense of \cite{CJ}?
\end{quest}

\begin{quest}
When $G$ is uniform, is every reflexive prime ideal of $\Omega_G$ principal?
\end{quest}

\subsection{The case when $G$ is almost simple} Recall that the compact $p-$adic analytic group $G$ is \emph{almost
simple} if every non-trivial closed normal subgroup of $G$ is open
(\ref{KdimOmega}). For such groups the constructions of
(\ref{ideals}) do not produce anything interesting because
$\Omega_G/w_H$ is artinian and hence finite dimensional over $\Fp$
for any closed normal subgroup $H \neq 1$. So Question G
specialises here to
\begin{quest}\label{AlmostSimple} Let $G$ be an almost simple uniform pro-$p$ group and let $P$ be a nonzero prime ideal of $\Omega_G$. Must $P$ be the unique maximal ideal of $\Omega_G$?
\end{quest}

We remind the reader that $x \in \Omega_G$ is \emph{normal} if $x \Omega_G = \Omega_G x$. Another closely related question is

\begin{quest} \label{NormalElts} Let $G$ be as in Question \ref{AlmostSimple}, with $G \ncong \Zp.$ Must any nonzero normal element of $\Omega_G$ be a unit?
\end{quest}

In \cite{Harris}, M. Harris claimed that, for $G$ as in Question
\ref{AlmostSimple}, any closed subgroup $H$ of $G$ with $2 \dim H
> \dim G$ gives rise to a non-zero two-sided ideal in $\Omega_G$,
namely the annihilator of the ``Verma module'' constructed by
induction from the simple $\Omega_H-$module. Unfortunately his
paper contains a gap, so Question \ref{AlmostSimple} remains open.
Some slight evidence towards a positive answer is provided by

\begin{thm}\cite[Theorem A]{Ard3} Suppose that $G$ is an almost simple uniform pro-$p$ group and that the Lie algebra of $G$ contains a copy of the two-dimensional non-abelian Lie algebra. Then for any two-sided ideal $I$ of $\Omega_G$,
\[\Kr(\Omega_G/I) \neq 1.\]
\end{thm}

Recall \cite[\S 6.4.4]{MR} that if $R$ is a oetherian ring with $\Kr(R)<\infty$, the \emph{classical Krull dimension} $\dim R$ of $R$ is the largest length of a chain of prime ideals of $R$. We always have $\dim R \leq \Kr(R)$; an easy consequence of the above result is
\[\dim (\Omega_G) < \dim G\]
whenever $G$ satisfies conditions of the Theorem.

\subsection{The case when $G$ is nilpotent} Towards the opposite end of the ``spectrum of
commutativity'' from the almost simple groups lie the nilpotent
groups. Motivated by analogous results for enveloping algebras of
nilpotent Lie algebras \cite[Chapter 4]{Dixm} and for group algebras $k[G]$ of finitely
generated nilpotent groups $G$ \cite[Theorem E]{Roseblade}, we ask

\begin{quest}\label{NilpControl} Let $G$ be a nilpotent uniform pro-$p$ group with
centre $Z$ and let $I$ be a nonzero ideal of $\Omega_G$. Does $I$
contain a non-zero central element? That is, is $I \cap \Omega_Z$
nonzero?
\end{quest}

S. J. Wadsley has shown that Question \ref{NilpControl} has an affirmative answer in the case when $G$ is one of the simplest possible nonabelian nilpotent uniform pro-$p$ groups:

\begin{thm}\cite[Theorem 4.10]{Wadsley} Let $G$ be a uniform Heisenberg pro-$p$ group with centre $Z$ and let $I$ be a nonzero two-sided ideal of $\Omega_G$. Then $I \cap \Omega_Z \neq 0$.
\end{thm}

A uniform pro-$p$ group $G$ is said to be \emph{Heisenberg} provided
its centre $Z$ is isomorphic to $\Zp$ and $G/Z$ is abelian. The main
idea of the proof of the above result is to show that for any integer
$t$, any finitely generated $\Omega_G-$module $M$ satisfying
$\Cdim(M) \leq \dim G/Z - t$ is actually finitely generated over
``most'' subalgebras $\Omega_H$ satisfying $Z \leq H$ and $\dim G/H = t$ \cite[Theorem 3.10]{Wadsley}.

In a more precise version of Question \ref{NilpControl}, one might
also hope that, when $G$ is nilpotent, ``small'' prime ideals $I$
in $\Omega_G$ are \emph{controlled} by $\Omega_Z$; that is
\[I = (I\cap \Omega_Z)\Omega_G.\] Question \ref{control} suggests
a more general version of this.

Moreover, one might even hope that arbitrary ideals of these
Iwasawa algebras of nipotent groups are constructed by means of a
sequence of centrally generated ideals - that is, one can ask:

\begin{quest}  Suppose that $G$ is a nilpotent uniform pro$-p$ group. If $I$ is an ideal of $\Omega_G$ strictly
 contained in $J(\Omega_G),$ is there a non-zero central element
 in $J(\Omega_G)/I$? \footnote{Compare with \cite[Proposition 4.7.1(i)]{Dixm}.}
 \end{quest}

\subsection{The case when $G$ is soluble}
Given the parallels pointed out in (\ref{skew}) between the
Iwasawa algebras of uniform soluble groups and the enveloping
algebras of finite dimensional complex soluble Lie algebras, it
is natural to wonder whether properties known for the latter case
might also be valid in the former. We give two sample questions of
this sort. Recall for the first that a prime ideal $P$ of the ring
$R$ is \emph{completely prime} if $R/P$ is a domain.

\begin{quest}
Let $G$ be a soluble uniform pro$-p$ group.
\begin{enumerate}[{(}i{)}]
\item Is every prime ideal of $\Omega_G$ completely prime?
\footnote{Compare with \cite[Theorem 3.7.2]{Dixm}.} \item Is the
prime spectrum of $\Omega_G$ the disjoint union of finitely many
commutative strata (along the lines of \cite[Theorem II.2.13]{BG},
but with necessarily non-affine strata)?
\end{enumerate}
\end{quest}

The simple possible nonabelian soluble case has been studied by O. Venjakob:

\begin{thm}\cite[Theorem 7.1]{Ven2} Let $G = X \rtimes Y$ be a nonabelian semidirect product of two copies of
$\Zp$. Then the only prime ideals of $\Omega_G$ are $0, w_X$
and $J(\Omega_G)$, and each one is completely prime. Moreover,
$w_X$ is generated by a normal element.
\end{thm}
An example of such a nonabelian semidirect product is provided by the group $B = \overline{\langle a_1,a_2 \rangle}$ considered in Example (\ref{SL2}).

Following J. E. Roseblade and D. S. Passman \cite[\S
1.5]{Roseblade}, we define the \emph{Zalesskii subgroup} $A$ of
the soluble uniform pro-$p$ group $G$ to be the centre of the
largest nilpotent closed normal subgroup $H$ of $G$. We say that
an ideal $I$ of $\Omega_G$ is \emph{faithful} if $G$ acts
faithfully on the quotient $\Omega_G/I$. If Question
\ref{NilpControl} has a positive answer, then it's possible that a
more general statement is true:

\begin{quest} \label{control} Let $G$ be a soluble uniform pro-$p$ group. Is every
faithful prime ideal of $\Omega_G$ controlled by the Zalesskii
subgroup $A$ of $G$?
\end{quest}


\begin{thebibliography}{99}
\bibitem{Ard} K. Ardakov, \emph{Krull dimension of Iwasawa Algebras}, J. Algebra \textbf{280} (2004), 190-206.
\bibitem{Ard2} K. Ardakov, \emph{The centre of completed group algebras of pro$-p$ groups}, Doc. Math \textbf{9} (2004), 599-606.
\bibitem{Ard3} K. Ardakov, \emph{Prime ideals in noncommutative Iwasawa algebras}, Math. Proc. Camb. Phil. Soc., to appear.
\bibitem{Ard4} K. Ardakov, \emph{Localisation at augmentation ideals in Iwasawa algebras}, submitted.
\bibitem{AB} K. Ardakov and K.A. Brown, \emph{Primeness, semiprimeness and localisation in Iwasawa algebras},
Transactions of the Amer. Math. Soc., to appear.
\bibitem{ASZ} K. Ajitabh, S. P. Smith, J. J. Zhang, \emph{Auslander-Gorenstein rings}, Comm. Algebra \textbf{26} (1998), 2159-2180.
\bibitem{BjE} J.-E. Bjork and E.K. Ekstrom, \emph{Filtered Auslander-Gorenstein rings}, in "Colloque en l'honneur de J. Dixmier," Birkhauser, Bosel, 1990.
\bibitem{Br} K. A. Brown, \emph{Height one primes of polycyclic group rings}, J. London Math. Soc. \textbf{32} (1985) 426-438; corrigendum J. London Math. Soc. .
\bibitem{BG} K. A. Brown and K.R. Goodearl, \emph{Lectures on Algebraic Quantum groups}, Birkhauser, 2002.
\bibitem{Bru} A. Brumer, \emph{Pseudocompact algebras, profinite groups and class formations}, J. Algebra \textbf{4} (1966), 442-470.
\bibitem{BH} W. Bruns and J. Herzog, \emph{Cohen-Macaulay Rings},
Cambridge Studies in Advanced Mathematics, CUP, (1993).
\bibitem{Carter} R. Carter, \emph{Simple groups of Lie type}, J. Wiley, London (1989).
\bibitem{Cau} G. Cauchon, \emph{Effacement des d\'erivations et spectres premiers des alg\`ebres quantiques},
 J. Algebra \textbf{260} (2003), no. 2, 476-518.
\bibitem{Ch} M. Chamarie, \emph{Modules sur les anneaux de Krull non commutatifs},
Paul Dubreil and Marie-Paule Malliavin algebra seminar, Lecture Notes in Math. vol. 1029, Springer, (1982),
283-310.
\bibitem{CJ} A.W. Chatters and D.A. Jordan, \emph{Non-commutative unique factorisation rings}, J. London Math. Soc. (2) 33 (1986), no. 1, 22--32.
\bibitem{C} J. Coates, Iwasawa algebras and arithmetic, S\'eminaire Bourbaki 2001/2002 Ast\'erisque 290 (2003), 37--52.
\bibitem{CSS} J. Coates, P. Schneider and R. Sujatha, \emph{Modules over Iwasawa algebras}, J. Inst. Math. Jussieu \textbf{2}, (2003) 73-108.
\bibitem{DDMS} J.D.Dixon, M.P.F. Du Sautoy, A.Mann, D.Segal, \emph{Analytic pro-p groups}, 2nd edition, CUP (1999).
\bibitem{Dixm} J. Dixmier, \emph{Enveloping Algebras}, Graduate Studies in Mathematics \textbf{11}, Amer. Math. Soc. (1996).
\bibitem{Harris} M. Harris, \emph{The annihilators of $p$-Adic Induced Modules}, J.Algebra \textbf{67}, 68-71 (1980).
\bibitem{L} M. Lazard, \emph{Groupes analytiques $p$-adiques}, Publ. Math. IHES \textbf{26} (1965), 389-603.
\bibitem{Le} T. Levasseur, \emph{Some properties of noncommutative regular graded rings}, \emph{Glasgow J. Math}, \textbf{34}, (1992) 277-300.
\bibitem{Le2} T. Levasseur, \emph{Krull dimension of the enveloping algebra of a semisimple Lie algebra},  Proc. Amer. Math. Soc.  \textbf{130}  (2002),  no. 12, 3519-3523.
\bibitem{LVO} L. Huishi and F. van Oystaeyen, \emph{Zariskian filtrations}, Kluwer Academic Publishers, K-monographs in Mathematics, vol. \textbf{2} (1996).
\bibitem{Martin} R. Martin, \emph{Skew group rings and maximal orders}, Glasgow Math. J. \textbf{37} (1995), no. 2, 249-263.
\bibitem{MaR} G. Maury and J. Raynaud, \emph{Ordres Maximaux au Sens de K. Asano} Lecture Notes in Math. vol. 808, Springer, 1980.
\bibitem{MR} J.C. McConnell, J.C. Robson, \emph{Noncommutative Noetherian rings}, AMS Graduate Studies
in Mathematics, vol. 30 (2001).
\bibitem{N} A. Neumann, \emph{Completed group algebras without zero divisors}, Arch. Math. (Basel) \textbf{51}, (1988) 496-499.
\bibitem{Pas} D.S. Passman, \emph{Infinite Crossed Products}, Pure and Applied Mathematics vol. \textbf{135}, Academic press, San Diego (1989).
\bibitem{Pas2} D.S. Passman, \emph{The Algebraic Structure of Group Rings}, Wiley-Interscience, New York, (1977).
\bibitem{Roseblade} J. E. Roseblade, \emph{Prime ideals in group rings of polycyclic groups}, Proc. London Math. Soc (3) \textbf{36}, (1978) 385-447.
\bibitem{Serre} J.-P. Serre, \emph{Sur la dimension homologique des groupes profinis}, Topology \textbf{3}, (1965) 413-420.
\bibitem{Smith} S. P. Smith, \emph{Krull dimension of factor rings of the enveloping algebra of a semisimple Lie algebra}, Math. Proc. Camb. Phil. Soc.  \textbf{93} (1983), no. 3, 459-466.
\bibitem{Staf} J. T. Stafford, \emph{Auslander-regular algebras and maximal orders},  J. London Math. Soc. (2) \textbf{50}  (1994),  no. 2, 276-292.
\bibitem{SZ} P. Samuel, O. Zariski, \emph{Commutative algebra}, Graduate Texts in Mathematics vol.\textbf{29}, Springer, 1975.
\bibitem{Ven1} O. Venjakob, \emph{On the structure theory of the Iwasawa algebra of a compact $p-$adic Lie group},  J. Eur. Math. Soc. (JEMS) 4 (2002), no. 3, 271-311.
\bibitem{Ven2}  O. Venjakob, \emph{A noncommutative Weierstrass Preparation Theorem and applications to Iwasawa Theory},  J. Reine Angew. Math. \textbf{559} (2003), 153-191.
\bibitem{VS} O. Venjakob and P. Schneider, \emph{On the dimension theory of skew power series rings}, preprint.
\bibitem{W} R. Walker, \emph{Local rings and normalising sets of elements}, Proc. London Math. Soc. \textbf{24} (1972), 27-45.
\bibitem{Wadsley} S. J. Wadsley, \emph{A Bernstein-type inequality for Heisenberg pro-$p$ groups}, preprint.
\bibitem{Wei} C. A. Weibel, \emph{An introduction to homological algebra}, Cambridge studies in advanced
mathematics \textbf{38}, CUP, (1994).
\bibitem{Zaks} A. Zaks, \emph{Injective dimension of semiprimary rings}, J. Algebra \textbf{13} (1969), 73-89.
\end{thebibliography}
\end{document}